\documentclass{amsart}
\usepackage{amsfonts,graphicx}

\makeatletter                                                  %
\def\th@plain{\slshape}                                        %
\makeatother                                                   %

\newcommand{\oi}{[0,1]}

\newcommand{\Nbb}{\mathbb{N}}
\newcommand{\Zbb}{\mathbb{Z}}
\newcommand{\Qbb}{\mathbb{Q}}
\newcommand{\Rbb}{\mathbb{R}}

\newcommand{\Tcal}{\mathcal{T}}

\newcommand{\Luk}{\L ukasiewicz}

\newcommand{\newword}[1]{\textsl{#1}}

\newcommand{\vect}[3]{#1_#2,\ldots ,#1_#3}

\newcommand{\angles}[1]{\langle #1 \rangle}

\DeclareMathSymbol{\upharpoonright}{\mathrel}{AMSa}{"16}
\let\restriction\upharpoonright
\DeclareMathSymbol{\nmid}{\mathrel}{AMSb}{"2D}

\DeclareMathOperator{\MaxSpec}{MaxSpec}
\DeclareMathOperator{\Cont}{Cont}

\DeclareMathOperator{\Free}{Free}

\DeclareMathOperator{\GL}{GL}
\DeclareMathOperator{\Mat}{Mat}

\theoremstyle{plain}
\newtheorem{theorem}{Theorem}[section]
\newtheorem{lemma}[theorem]{Lemma}

\theoremstyle{definition}

\begin{document}

\bibliographystyle{plain}

\sloppy

\title[Bernoulli automorphisms]{Bernoulli automorphisms\\
of finitely generated free MV-algebras}

\author[G. Panti]{Giovanni Panti}
\address{Department of Mathematics and Computer Science\\
University of Udine\\
Via delle Scienze 208\\
33100 Udine, Italy}
\email{panti@dimi.uniud.it}

\begin{abstract}
MV-algebras can be viewed either as the Lindenbaum algebras of \Luk\ infinite-valued logic, or as unit intervals $[0,u]$ of lattice-ordered abelian groups in which a strong order unit $u>0$ has been fixed.
They form an equational class, and the free $n$-generated free MV-algebra is representable as an algebra of piecewise-linear continuous functions with integer coefficients over the unit $n$-dimensional cube.
In this paper we show that
the automorphism group of such a free algebra contains elements having strongly chaotic behaviour, is the sense that their duals are measure-theoretically isomorphic to a Bernoulli shift. This fact is noteworthy from the viewpoint of algebraic logic, since it gives a distinguished status to Lebesgue measure as an averaging measure on the space of valuations. As an ergodic theory fact, it provides explicit examples of volume-preserving homeomorphisms of the unit cube which are piecewise-linear with integer coefficients, preserve the denominators of rational points, and enjoy the Bernoulli property.
\end{abstract}


\thanks{\emph{2000 Math.~Subj.~Class.}: 06D35; 37A05}

\maketitle

\section{Preliminaries}

An \newword{MV-algebra} is an algebra $(A,\oplus,\neg,0)$ such that $(A,\oplus,0)$ is a commutative monoid and the identities $\neg\neg f=f$, $f\oplus\neg 0=\neg 0$, and $\neg(\neg f\oplus g)\oplus g=\neg(\neg g\oplus f)\oplus f$ are satisfied. MV-algebras can be viewed either as the Lindenbaum algebras of \Luk\ infinite-valued logic~\cite{chang1}, \cite{chang2}, or as unit intervals $[0,u]$ of lattice-ordered abelian groups in which a strong order unit $u>0$ has been fixed~\cite{mundicijfa}.
We assume some familiarity with the basics of the theory; the book~\cite{CignoliOttavianoMundici00} provides an extensive treatment, but~\cite{mundici95} or~\cite[\S2-4]{mundicijfa} cover material sufficient for our needs.

We are concerned with automorphisms of the free MV-algebra over $n$ generators, which we denote by $\Free_n$. Such an algebra has a representation as a space of piecewise-linear 
continuous functions with integer coefficients, as follows.
A \newword{rational cellular complex over the $n$-cube} is a finite set $W$ of cells (i.e., compact convex polyhedrons) whose union is $\oi^n$ and such that:
\begin{enumerate}
\item every vertex of every cell of $W$ has rational coordinates;
\item if $C\in W$ and $D$ is a face of $C$, then $D\in W$;
\item every two cells intersect in a common face.
\end{enumerate}

Equip the real unit interval $\oi$ with the operations $a\oplus b=\min(1,a+b)$, $\neg a=1-a$, and the constant $0$. Then $\oi$ is an MV-algebra, and so is the power $\oi^{\oi^n}$ under componentwise operations. By mapping the $i$th free generator to the $i$th projection function $x_i:\oi^n\to\oi$, we obtain a representation $\rho:\Free_n\to\oi^{\oi^n}$, which is faithful by Chang's completeness theorem~\cite{chang2}, \cite{pantilu}. McNaughton's theorem~\cite{mcnaughton}, \cite{mundicimn} states then that the range of $\rho$ is exactly the set of those continuous functions $f:\oi^n\to\oi$ for which there exists a complex 
as above and affine linear functions $F_j(\bar x)=a_j^1x_1+\cdots+a_j^nx_n+a_j^{n+1}\in\Zbb[\vect x1n]$, in 1-1 correspondence with the $n$-dimensional cells $C_j$ of the complex, such that $f\restriction C_j=F_j$ for each $j$. In the following we will tacitly identify $\Free_n$ with its image under $\rho$, except when discussing states in \S2.

Let $\sigma$ be an endomorphism of $\Free_n$. The \newword{dual} of $\sigma$ is the mapping $S:\oi^n\to\oi^n$ defined by $S(p)_i=(\sigma(x_i))(p)$. 
This is equivalent to saying that there exists a complex $W$ as above and matrices
$P_j\in\Mat_{n+1}(\Zbb)$, in 1-1 correspondence with the $n$-dimensional cells $C_j$ of $W$, such that:
\begin{enumerate}
\item every $P_j$ has last row $(0\ldots0\,1)$;
\item $P_j$ expresses $S\restriction C_j$ in homogeneous coordinates (i.e., if $p=(\vect r1n)\in C_j$ and
$S(p)=(\vect s1n)$, then $(s_1\ldots s_n\,1)^{tr}=
P_j(r_1\ldots r_n\,1)^{tr}$).
\end{enumerate}
As in the proof of~\cite[Theorem~3.4]{mundici95}, $\sigma$ is an automorphism iff $S$ is injective and all $P_j$'s are in $\GL_{n+1}(\Zbb)$. In particular, all duals of automorphisms preserve Lebesgue measure $\lambda$ on $\oi^n$; we call these duals \newword{McNaughton homeomorphisms} of the unit cube. We say that an automorphism is \newword{Bernoulli} if the corresponding McNaughton homeomorphism is isomorphic to a two-sided Bernoulli shift. This means that there exists a probability space $(X,\mu)$ and a measure-theoretic isomorphism $\varphi:(\oi^n,\lambda)\to(X^\Zbb,\mu^\Zbb)$ such that $\varphi\circ S=B\circ\varphi$, where $B$ is the shift on~$X^\Zbb$.

\section{Statement of the result}

The main result of this paper is the construction of an example (actually, a family of examples) of Bernoulli automorphisms of $\Free_2$. This is the content of Theorem~\ref{ref1}, whose statement and proof use some amount of ergodic theory, for which we refer to~\cite{CornfeldFomSi82} or~\cite{Walters82}. In this section we explain the significance of our result both from the viewpoint of ergodic theory and from that of algebraic logic.

\begin{theorem}
There\label{ref1} is an explicitly constructible family $\{\beta_{lm}:1\le l,m\in\Nbb\}$ of automorphisms of $\Free_2$ such that, for every $\beta_{lm}$ in the family, the dual McNaughton homeomorphism $B_{lm}$ has the following properties:
\begin{enumerate}
\item $B_{lm}$ is ergodic with respect to Lebesgue measure on the $2$-cube;
\item $B_{lm}$ is isomorphic to a Bernoulli shift;
\item $B_{lm}$ fixes pointwise the boundary of $\oi^2$ and permutes the set of points having rational coordinates and given denominator;
\item $B_{lm}$ is chaotic in the sense of Devaney (i.e., is topologically transitive, has sensitive dependence on initial conditions, and the set of points having finite orbit is dense in $\oi^2$).
\end{enumerate}
\end{theorem}

It is a classic result of Oxtoby and Ulam that in the space of measure-preserving homeomorphisms of the unit cube with the uniform topology the subset of ergodic ones is $G_\delta$ comeager. The same property holds for the subset of chaotic homeomorphisms; see~\cite{AlpernPrasad00} for a presentation of these results. In particular, we are assured of the existence of many (from the topological point of view) such homeomorphisms. Note however that these theorems are obtained by Baire category methods, and do not provide us with any explicit example. To the best of our knowledge, in the literature one finds two main examples:
\begin{enumerate}
\item The topological skew products in~\cite{Xu90}, obtained by conjugating the Besicovitch mapping. They are topologically transitive homeomorphisms that do not preserve Lebesgue measure and restrict to irrational rotations on the boundary.
\item Quotients of Anosov diffeomorphisms of the $2$-torus give chaotic homeomorphisms of the sphere, which in turn yield chaotic homeomorphisms of $\oi^2$ by blowing up a fixed point. The resulting maps have two fixed points on the boundary (one attracting and one repelling, corresponding to the eigenspaces of the original diffeomorphism). This construction is modified in~\cite{Katok79}, obtaining Bernoulli diffeomorphisms of the unit disk. Another modification is in~\cite{Cairns_et_al95}, where the authors start from a linked twist of the torus, and obtain a chaotic homeomorphism of $\oi^2$ which fixes the boundary.
\end{enumerate}

Theorem~\ref{ref1} allows us to add our family $B_{lm}$ to the above list: it is worth noting that the arithmetical nature of McNaughton homeomorphisms ---namely, their piecewise-linearity, and the restriction to integer coefficients--- makes them much more rigid objects than the examples (1) and (2), which exploit the full power of the real numbers.

From the viewpoint of algebraic logic, our result gives a distinguished status to Lebesgue measure as an averaging measure on the space $\oi^n$ of truth-value assignments. Recall, first of all, the representation $\rho$ described in the previous section. Let $f(\vect x1n)$ be an $n$-variable formula in \Luk\ infinite-valued logic, which we identify with an element of $\Free_n$. We assign to $f$ an ``average truth-value'' $m(f)$ by setting
\begin{equation}\tag{A}
m(f)=\int_{\oi^n}\rho(f)\,d\lambda,
\end{equation}
where $\lambda$ is Lebesgue measure on the $n$-cube. We have $m(0)=0$, $m(1)=1$, and $m(f\oplus g)=m(f)+m(g)$ provided that $\neg(\neg f\oplus\neg g)=0$: hence the map $m:\Free_n\to\oi$ is an example of a \newword{state}, as defined in~\cite{mundici95}. By Theorem~\ref{ref53} below $m$ does not depend on the choice of a free generating set $\{\vect x1n\}$ for $\Free_n$.

\begin{theorem}
Let $n\ge2$.
There\label{ref53} exists exactly one state on $\Free_n$, namely $m$, that satisfies the following two conditions:
\begin{itemize}
\item[(B1)] $m(\sigma(f))=m(f)$, for every automorphism $\sigma$ of $\Free_n$;
\item[(B2)] if, for every $g\in\Free_n$, $f\land g=0$ implies $g=0$, then
$$
\lim_{k\to\infty} m(\underbrace{f\oplus\cdots\oplus f}_\text{$k$ summands})=1.
$$
\end{itemize}
\end{theorem}
\begin{proof} By~\cite[Proposition~1.1]{pantiinvariant}, for every MV-algebra $A$ the set of states on $A$ is in 1--1 correspondence with the set of regular Borel probability measures on $\MaxSpec A$, the latter being the space of maximal ideals of $A$ equipped with the hull-kernel topology. Write $\Cont X$ for the MV-algebra of all $\oi$-valued continuous functions on the topological space $X$, under pointwise operations. The dual of the faithful representation $\rho:\Free_n\to\Cont\oi^n$ is the map
\begin{equation*}
\begin{split}
R:\oi^n &\to \MaxSpec\Free_n \\
t &\mapsto \{f\in\Free_n:\bigl[\rho(f)\bigr](t)=0\},
\end{split}
\end{equation*}
which is a homeomorphism by~\cite[Lemma~8.1]{mundicijfa}. It follows that the states on $\Free_n$ are in 1--1 correspondence with the 
Borel probability measures on $\oi^n$ (which are necessarily regular); in this bijection, $m$ corresponds to Lebesgue measure $\lambda$, by~(A). Hence, our claim boils down to the fact that Lebesgue measure is the only Borel measure $\mu$ on $\oi^n$ satisfying the analogues of~(B1) and~(B2) (in which, of course, $\Free_n$ is replaced by $\rho[\Free_n]$ and $m$ by integration w.r.t.~$\mu$). This fact is proved in~\cite[Theorem~2.4]{pantiinvariant}.
\end{proof}

The restriction $n\ge2$ in Theorem~\ref{ref53} is irrelevant, since we can always add a dummy variable to any 1-variable formula; see the discussion about coherence in~\cite[\S 2]{pantiinvariant}. On the other hand, it is crucial to note that, although the proof of \cite[Theorem~2.4]{pantiinvariant} exploits a fair amount of the piecewise-linear structure of $\oi^n$, Theorem~\ref{ref53} is an intrinsic, purely algebraic characterization of $m$. 
As a matter of fact, we \emph{define} $m$ as the only state on $\Free_n$ that satisfies~(B1) and~(B2).
Note also that, given any representation of $\Free_n$ as an MV-algebra of continuous functions over a topological space $X$, Theorem~\ref{ref53} uniquely determines a distinguished probability measure on $X$, regardless of any arithmetical or piecewise-linear structure $X$ may or may not have.

Let us further clarify this point: let $X$ be any compact Hausdorff space, and let $\theta$ be any faithful representation of $\Free_n$ over $X$ (i.e., $\theta$ is an injective homomorphism of $\Free_n$ onto a subalgebra of $\Cont X$). Without loss of generality, $\theta[\Free_n]$ separates the points of $X$ (if not, replace $X$ by an appropriate quotient space). By~\cite[Theorem~3.4.3]{CignoliOttavianoMundici00}, the map
\begin{equation*}
\begin{split}
T:X &\to \MaxSpec\Free_n \\
x &\mapsto \{f\in\Free_n:\bigl[\theta(f)\bigr](x)=0\}
\end{split}
\end{equation*}
is a continuous bijection.
The state $m$ determines then a measure $\mu$ on $X$, namely the pull-back along $T$ of the measure on $\MaxSpec\Free_n$ associated to $m$. Note that $\mu$ is Borel regular with respect to the topology $\Tcal$ pulled back by $T$ ($\Tcal$ may be coarser than the original topology on $X$). In equivalent terms, $\mu$ can either be viewed as the pull-back along $R^{-1}\circ T$ of Lebesgue measure on $\oi^n$, or as the unique $\Tcal$-Borel regular measure on $X$ that satisfies the following:
$$
\int_X\theta(f)\,d\mu=\int_{\oi^n}\rho(f)\,d\lambda,
\qquad\text{for every $f\in\Free_n$.}
$$
As a simple example, let $n=1$, $X=\oi$, and let $\theta$ map the free generator $x_1$ of $\Free_1$ to the function $\bigl[\theta(x_1)\bigr](x)=x^3$. We let $x$ and $t$ vary in $X$ and in $\oi$, respectively; of course $\bigl[\rho(x_1)\bigr](t)=t$.
By standard duality, $R^{-1}\circ T$ must obey
$$
\bigl[\theta(x_1)\bigr](x)=
\bigl[\rho(x_1)\bigr]\bigl((R^{-1}\circ T)(x)\bigr),
$$
and hence $(R^{-1}\circ T)(x)=x^3$. It follows that the $\mu$-measure of, say, the interval $[\sqrt[3]{t_1},\sqrt[3]{t_2}]\subseteq X$ is the $\lambda$-measure of the interval $[t_1,t_2]\subseteq\oi$. In this specific example $R^{-1}\circ T$ is differentiable, and we can recast things in terms of volume forms. To compute $m(x_1)$ we integrate $\bigl[\theta(x_1)\bigr](x)$ not against the form $dx$, but against the image of $dt$ under the cotangent map induced by $R^{-1}\circ T$, namely $3x^2\,dx$, and we obtain, of course, the value $1/2$.

Having completed our discussion about the state $m$, we resume our standing identification of $\Free_n$ with $\rho[\Free_n]$. 
Let $q$ be any density on $\oi^n$ (i.e., any measurable function from the unit cube to the nonnegative real numbers whose integral w.r.t.~$\lambda$ is $1$). Then $q$ induces a measure $\mu_q$ on $\oi^n$ by $d\mu_q=q\,d\lambda$, and the function $m_q$ assigning to $f$ the value of the integral of $\rho(f)$ w.r.t.~$\mu_q$ is again a state. By possibly adding a dummy variable, we may assume that $n$ is even; since the Bernoulli property is preserved under direct products~\cite[Ch.~10 \S1]{CornfeldFomSi82}, Theorem~\ref{ref1} immediately guarantees the existence of a Bernoulli automorphisms $\sigma$ of $\Free_n$.
The dual $S$ of $\sigma$ is mixing~\cite[Corollary~4.33.1]{Walters82}, and hence the push forward $S_*^k\mu_q$ of $\mu_q$ by $S^k$ converges to $\lambda$ in the weak${}^*$ topology~\cite[\S4.9 and Theorem~6.12(ii)]{Walters82}. In particular we get
\begin{multline*}
\lim_{k\to\infty}m_q(\sigma^k(f))=
\lim_{k\to\infty}\int\rho(\sigma^k(f))\,d\mu_q=
\lim_{k\to\infty}\int \rho(f)\circ S^k\,d\mu_q\\
=\lim_{k\to\infty}\int \rho(f)\,dS_*^k\mu_q=
\int \rho(f)\,d\lambda=
m(f).
\end{multline*}
Hence the action of our family makes $m$ an attractor for the family of all states induced by a density: every ``wrong'' average truth-value assignment eventually converges to the right one.


\section{The family $B_{lm}$}

Consider the following two complexes over $\oi^2$:
\begin{figure}[!h]
\begin{center}
\includegraphics[height=4cm,width=4cm]{figura.1}
\qquad
\includegraphics[height=4cm,width=4cm]{figura.2}
\end{center}
\end{figure}

In the first complex, the vertices of the lower inner triangle are $p_0=(1/4,1/4)$, $p_1=(1/2,1/4)$, $p_2=(1/4,1/2)$, and the picture is symmetric under a $\pi$ rotation about $(1/2,1/2)$. For $0\le i\le 2$, let $p_i'$ be the vertex symmetric to $p_i$. Then there exists a unique homeomorphism $R_1$ such that:
\begin{enumerate}
\item $p_i$ is mapped to $p_{i+1\pmod{3}}$, and analogously for $p_i'$;
\item every other vertex is fixed;
\item $R_1$ is affine linear on each cell.
\end{enumerate}
By direct computation one sees that $R_1$ is a McNaughton homeomorphism.

In the second complex, the innermost square has vertices $q_0=(3/8,3/8)$, $q_1=(5/8,3/8)$, $q_2=(5/8,5/8)$, $q_3=(3/8,5/8)$, while the intermediate one has vertices $p_0$, $p_3=(3/4,1/4)$, $p_0'$, $p_4=(1/4,3/4)$ (see the picture of the intermediate square before the proof of Lemma~\ref{ref7}). Define $R_2$ to be the unique homeomorphism that maps $q_i$ to $q_{i-1\pmod4}$, fixes every other vertex, and is affine linear on each cell. $R_2$~is not a McNaughton homeomorphism since, e.g., the matrix that maps the triangle $\langle p_0,q_1,q_0\rangle$ to $\langle p_0,q_0,q_3\rangle$ has rational non-integer entries.

\begin{lemma}
$R_2^2$ is a McNaughton homeomorphism.
\end{lemma}
\begin{proof}
Let $B$, $C$ be $2$-dimensional cells in the second complex, with corresponding matrices $P$, $Q$, respectively. It suffices to show that if the topological interior of $B\cap R_2^{-1}C$ is not empty, then $QP$ has integer entries. Observe that the complex is symmetric under a clockwise $\pi/2$ rotation about $(1/2,1/2)$, and this rotation is expressible in homogeneous coordinates via a matrix $A\in\GL_3(\Zbb)$. We can clearly restrict our attention to the triangles in the ring bounded by the innermost and the intermediate squares. We then see easily that, fixing $B=\langle p_0,p_3,q_1\rangle$ and $C=\langle p_0,q_1,q_0\rangle$, we only have to check that the four matrices $P^2,QP,(APA^{-1})Q,(AQA^{-1})Q$ have integer entries, and this can be done by hand computation.
\end{proof}

For every $1\le l,m$, let $B_{lm}=R_2^{4m}\circ R_1^{3l}$, and let $\beta_{lm}$ be the associated automorphism of $\Free_2$. The family of all $\beta_{lm}$'s is the one in the statement of Theorem~\ref{ref1}.

\section{Proof of Theorem~\ref{ref1}}

We fix $l,m\ge1$, and let $\beta=\beta_{lm}$, $B=B_{lm}$.
We first observe that $\beta$ is explicitly constructible. Indeed, by using the theory of Schauder hats~\cite[p.~598]{mundicimn}, \cite[\S3]{pantilu}, we can write down algorithmically syntactic expressions for the formulas $t_1=x_1\circ R_1$, $t_2=x_2\circ R_1$, $t_3=x_1\circ R_2^2$, $t_4=x_2\circ R_2^2$. The images $\beta(x_1)$ and $\beta(x_2)$ of the free generators $x_1,x_2$ of $\Free_2$ are then (huge) formulas obtained from $\vect t14$ via recursive nesting.

Since $R_1$ and $R_2$ fix the boundary of $\oi^2$, so does $B$.
Define the \newword{denominator} of the rational point $p=(r_1,r_2)\in\oi^2\cap\Qbb^2$ to be the unique integer $d\ge1$ such that $dr_1,dr_2,d$ are relatively prime integers.
From the definitions, it is immediate that every McNaughton homeomorphism permutes the set $A_d$ of all points having denominator $d$. Since every $A_d$ is finite, and the union of all $A_d$'s is the set of rational points, which is dense in the square, it follows that every topologically transitive McNaughton homeomorphism is chaotic in the sense of Devaney~\cite[\S4.5]{AlpernPrasad00}; the statements~(3)
and~(4) in Theorem~\ref{ref1} are then consequences of~(1) and~(2). 
Note that, since McNaughton homeomorphisms are determined by their actions on the $A_d$'s, the automorphism group of $\Free_2$ is residually finite, and hence acts chaotically on some Hausdorff space~\cite[Theorem~1]{Cairns_et_al95}. Theorem~\ref{ref1}(4) says that in this group there are specific elements that act chaotically, and that the space on which the action takes place may be assumed to be the unit square.

To prove the statements~(1) and~(2) in Theorem~\ref{ref1}, we use Pesin theory~\cite{pesin77}: the key ingredients will be Lemma~\ref{ref2} and Lemma~\ref{ref3}, stated below.

\begin{lemma}
There\label{ref2} is a $B$-invariant and eventually strictly invariant cone field $U_q$, defined $\lambda$-almost everywhere on $\oi^2$, and depending measurably on $q$.
\end{lemma}

By~\cite[p.~152]{wojtkowski85}, \cite[Theorem~III.4.4]{chernovmar01}, Lemma~\ref{ref2} entails that for $\lambda$-all points $q\in\oi^2$ the tangent space at $q$ decomposes in the direct sum $E_q^s\oplus E_q^u$ of stable and unstable subspaces, the Lyapunov exponents are defined at $q$, and their sum is $0$. These facts determine foliations of the square into local stable and unstable curves $\gamma_q^s$ and $\gamma_q^u$. Due to the piecewise linearity of $B$, these (un)stable leafs are indeed segments, as in~\cite[Proposition~1]{wojtkowski80}, and the resulting foliations are absolutely continuous~\cite[p.~71]{wojtkowski80}. Using~\cite[Theorems~7.2, 7.9, 8.1]{pesin77} and the argument in~\cite{przytycki83}, Theorem~\ref{ref1} will result from the following lemma.

\begin{lemma}
For\label{ref3} $\lambda$-all $p,q\in\oi^2$, there exist $h_0,k_0\ge0$ such that, for every $h\ge h_0$ and $k\ge k_0$, we have
$$
B^{-h}[\gamma_p^s]\cap B^k[\gamma_q^u]\not=\emptyset.
$$
\end{lemma}

We prove the above lemmas in the next two sections.

\section{Proof of Lemma~\ref{ref2}}

Let $D$ be the unit disk in polar coordinates $(r,\theta)$. Let $f:\oi\to\oi$ be the continuous piecewise-fractional function defined by
$$
f(r)=\min\bigl(1,(1/r-1)/3\bigr),
$$
and let $F:D\to D$ be the twist mapping $F(r,\theta)=(r,\theta+2\pi lf(r))$.
Let $L\subset\oi^2$ be the triangle with vertices $(0,0), (1,0), (0,1)$, and let $M:D\to L$ be the homeomorphism defined by
\begin{equation*}
\begin{split}
M(0,\theta) &= (1/3,1/3); \\
M(1,\theta)  &= \begin{cases}
\bigl(1-3\theta/2\pi,3\theta/2\pi\bigr), & 
\text{if } 0\le\theta<2\pi/3; \\
\bigl(0,2-3\theta/2\pi\bigr), & 
\text{if } 2\pi/3\le\theta<4\pi/3; \\
\bigl(3\theta/2\pi-2,0\bigr), & 
\text{if } 4\pi/3\le\theta<2\pi;  
\end{cases} \\
M(r,\theta) &= (1-r)M(0,\theta)+rM(1,\theta).
\end{split}
\end{equation*}

\begin{lemma}
$M$\label{ref4} conjugates $F$ with $R_1^{3l}\restriction L$ (i.e., $M\circ 
F=R_1^{3l}\circ M$).
\end{lemma}
\begin{proof}
This is a straightforward computation, which we just sketch. We fix without loss of generality $1/4<r<1$. By definition of $R_1$, the points of the triangle $T(r)=\{M(r,\theta):0\le\theta<2\pi\}$ ``rotate counterclockwise'' along $T(r)$ under the action of $R_1$. We have to check that the speed of the rotation is $f(r)/3$, upon normalizing the lengths of the sides of $T(r)$ to be $r/3$. Since $R_1$ is piecewise affine, we may do the checking only for the $6$ points in the intersection of $T(r)$ with the $1$-dimensional simplices in the complex defining $R_1$.
Consider, e.g., the point $w$ in $T(r)\cap\angles{(0,0),p_2}$. It moves to $z$ in $T(r)\cap\angles{(0,0),p_0}$ with speed $\overline{wz}/(\text{length of $T(r)$})=
\overline{wz}/r$. Since $\overline{wz}/(1-r)=\overline{p_2p_0}\big/(3/4)=(4^{-1}/3)\big/(3/4)=1/9$ (because the triangles $\angles{(0,0),w,z}$ and $\angles{(0,0),p_2,p_0}$ are similar, and $p_0,p_2\in T(4^{-1})$), we obtain $\overline{wz}/r=(1-r)/(9r)=f(r)/3$, as desired.
\end{proof}

If $A$ is a $\pi$ rotation of the unit square about $(1/2,1/2)$, $L'=AL$, and $M'=A\circ M$, then clearly $M'$ conjugates $F$ with $R_1^{3l}\restriction L'$.

Define now $g:\oi\to\oi$ and $G:D\to D$ by
\begin{equation*}
\begin{split}
g(r) &= \min(1,1/r-1), \\
G(r,\theta) &= (r,\theta-2\pi mg(r)).
\end{split}
\end{equation*}
Let $E\subset\oi^2$ be the square $\langle p_0,p_3,p_0',p_4\rangle$, and let $N:D\to E$ be the homeomorphism
$$
N(r,\theta)=(1/2,1/2)+\begin{cases}
r(1/4, -1/4+\theta/\pi), & \text{if $0\le\theta<\pi/2$;} \\
r(3/4-\theta/\pi,1/4), & \text{if $\pi/2\le\theta<\pi$;} \\
r(-1/4,5/4-\theta/\pi), & \text{if $\pi\le\theta<3\pi/2$;} \\
r(-7/4+\theta/\pi,-1/4), & \text{if $3\pi/2\le\theta<2\pi$.}
\end{cases}
$$
The mapping $N$ arises from the same construction as the mapping $M$: we have just given directly its analytic expression, as the reader can easily verify. The same proof as in Lemma~\ref{ref4} shows that $N$ conjugates $G$ with $R_2^4\restriction E$.

Consider the annuli $H=\{(r,\theta):1/4<r<1\}$, $K=\{(r,\theta):1/2<r<1\}$ and the standard row vector basis 
$(\partial/\partial r\rvert_{(r,\theta)} \;
\partial/\partial \theta\rvert_{(r,\theta)})$
on their tangent spaces. The differentials $T_{(r,\theta)}F$ and $T_{(r,\theta)}G$ have then the matrix form
\begin{align*}
T_{(r,\theta)}F &= \begin{pmatrix}
1 & 0 \\
-2\pi l/3r^2 & 1
\end{pmatrix}, &
T_{(r,\theta)}G &= \begin{pmatrix}
1 & 0 \\
2\pi m/r^2 & 1
\end{pmatrix}.
\end{align*}

We introduce the matrices
\begin{align*}
A_1(r) &= \begin{pmatrix}
1 & 0 \\
\pi l/3r^2 & -1
\end{pmatrix}, &
A_2(r) &= \begin{pmatrix}
1 & 0 \\
-\pi l/3r^2 & -1
\end{pmatrix}, \\
A_3(r) &= \begin{pmatrix}
1 & 0 \\
-\pi m/r^2 & 1
\end{pmatrix}, &
A_4(r) &= \begin{pmatrix}
1 & 0 \\
\pi m/r^2 & 1
\end{pmatrix}.
\end{align*}
For each $i=1,2$ (respectively, $i=3,4$), let $U^i_{(r,\theta)}$ be the cone in $T_{(r,\theta)}H$ (respectively, in $T_{(r,\theta)}K$)
of all vectors whose coordinates in basis $(\partial/\partial r\vert_{(r,\theta)}\quad
\partial/\partial\theta\vert_{(r,\theta)})A_i(r)$ are both $\ge0$ or both $\le0$. Then $T_{(r,\theta)}F$ maps bijectively $U^1_{(r,\theta)}$ onto $U^2_{F(r,\theta)}$, and $T_{(r,\theta)}G$ maps bijectively $U^3_{(r,\theta)}$ onto $U^4_{G(r,\theta)}$.
In particular, since $U^2_{(r,\theta)}\subset U^1_{(r,\theta)}$ and
$U^4_{(r,\theta)}\subset U^3_{(r,\theta)}$, the differential $T_{(r,\theta)}F$ maps $U^2_{(r,\theta)}$ into $U^2_{F(r,\theta)}$, and $T_{(r,\theta)}G$ maps $U^4_{(r,\theta)}$ into $U^4_{G(r,\theta)}$. Recall that a cone is \newword{strictly contained} in another cone if the former is contained ---except for the common vertex--- in the topological interior of the latter.
Consider the parallelograms $\vect P14$ in the following picture:
\begin{figure}[!h]
\begin{center}
\includegraphics[height=5cm,width=5cm]{figura.3}
\end{center}
\end{figure}

\begin{lemma}
Let\label{ref7} $q=(x,y)$ be in the topological interior of $P_1\cup P_2$. Then the cone $(T_{M^{-1}q}M)[U^2_{M^{-1}q}]$ is strictly contained in 
$(T_{N^{-1}q}N)[U^3_{N^{-1}q}]$, and $(T_{N^{-1}q}N)[U^4_{N^{-1}q}]$ is strictly contained in 
$(T_{M^{-1}q}M)[U^1_{M^{-1}q}]$.
\end{lemma}
\begin{proof}
Four computations are needed, according whether $q$ is in $P_1$ or $P_2$: it takes about one hour to perform all four using a computer algebra system. We will carry on just one computation, the others being analogous. Namely, we will show that
$(T_{N^{-1}q}N)[U^4_{N^{-1}q}]$ is strictly contained in 
$(T_{M^{-1}q}M)[U^1_{M^{-1}q}]$, on the assumption that
$q=M(r_1,\theta_1)=N(r_2,\theta_2)$ is in $P_1$; hence $0<\theta_1<2\pi/3$ and $3\pi/2<\theta_2<2\pi$
(see the picture above, in which the cones are computed for $l=m=1$).
Introducing the basis 
$(\partial/\partial x\rvert_{(x,y)} \;
\partial/\partial y\rvert_{(x,y)})$
on $T\oi^2$,
the differentials $T_{(r_1,\theta_1)}M$ and $T_{(r_2,\theta_2)}N$ are given by the matrices
\begin{align*}
B_1 &= \begin{pmatrix}
2/3 - 3\theta_1/2\pi & -3r_1/2\pi \\
-1/3+3\theta_1/2\pi & 3r_1/2\pi
\end{pmatrix}, &
B_2 &= \begin{pmatrix}
-7/4+\theta_2/\pi & r_2/\pi \\
-1/4 & 0
\end{pmatrix}.
\end{align*}

We then have to check that the images of the columns of $A_4(r_2)$ under $B_2$ have coordinates all $>0$ or all $<0$ in the basis given by the images of the columns of $A_1(r_1)$ under $B_1$: in other words, that the matrix
$$
\bigl[B_1A_1(r_1)\bigl]^{-1}\bigl[B_2A_4(r_2)\bigr]
$$
has all entries $>0$ or all $<0$. Observe that
$$
A_1(r_1)^{-1}=\begin{pmatrix}
1 & 0 \\
a & 1
\end{pmatrix}
\begin{pmatrix}
1 & 0 \\
1/r_1^2 & -1
\end{pmatrix}
$$
for a certain $a>0$ ($a$ depends on $r_1$ and $l$), and analogously
$$
A_4(r_2)=\begin{pmatrix}
1 & 0 \\
3/r_2^2 & 1
\end{pmatrix}
\begin{pmatrix}
1 & 0 \\
b & 1
\end{pmatrix}
$$
for some $b>0$, depending on $r_2$ and $m$. Hence, it is sufficient to check the entries of
$$
\begin{pmatrix}
1 & 0 \\
1/r_1^2 & -1
\end{pmatrix}
B_1^{-1}B_2
\begin{pmatrix}
1 & 0 \\
3/r_2^2 & 1
\end{pmatrix}.
$$

We can multiply our matrices by factors $>0$ to get rid of denominators: in particular, we compute
\begin{equation}\tag{$*$}
3r_1^2\begin{pmatrix}
1 & 0 \\
1/r_1^2 & -1
\end{pmatrix}
B_1^{-1}=
\begin{pmatrix}
9r_1^2 & 9r_1^2 \\
9r_1\theta_1-2r_1\pi+9 & 9r_1\theta_1-4r_1\pi+9
\end{pmatrix}.
\end{equation}
From the definition of $M$ we have the identity
$$
\begin{pmatrix}
x-1/3 \\
y-1/3
\end{pmatrix}
=\begin{pmatrix}
2/3 & -3/2\pi \\
-1/3 & 3/2\pi
\end{pmatrix}
\begin{pmatrix}
r_1 \\
r_1\theta_1
\end{pmatrix},
$$
which we invert and substitute in the right-hand side of~($*$), obtaining
$$
C=\begin{pmatrix}
\begin{matrix}
81x^2+81y^2+162xy-\\
-108x-108y+36
\end{matrix}
& \begin{matrix}
81x^2+81y^2+162xy-\\
-108x-108y+36
\end{matrix} \\
\noalign{\bigskip} 
6\pi y-2\pi+9 & -6\pi x+2\pi+9
\end{pmatrix}.
$$
Analogously, we compute
$$
2\pi r_2B_2\begin{pmatrix}
1 & 0 \\
3/r_2^2 & 1
\end{pmatrix}
$$
and apply to the result the substitution $r_2=-4y+2$, $r_2\theta_2=\pi x-7\pi y+3\pi$ that follows from inverting the identity
$$
\begin{pmatrix}
x-1/2 \\
y-1/2
\end{pmatrix}
=\begin{pmatrix}
-7/4 & 1/\pi \\
-1/4 & 0
\end{pmatrix}
\begin{pmatrix}
r_2 \\
r_2\theta_2
\end{pmatrix}.
$$
We call $D$ the resulting matrix, namely
$$
D=\begin{pmatrix}
2\pi x-\pi+6 & 32y^2-32y+8 \\
2\pi y-\pi & 0
\end{pmatrix}.
$$
We finally compute $CD$ and apply to the result the change of variables $x=z-w+1/2$, $y=w+1/4$, so that $q=(z,w)$, with $0<z<1/4$ and $0<w<1/8$. We thus get the matrix
$$
\begin{pmatrix}
\begin{matrix}
(9/32)(12z+1)^2(4\pi z-\pi+12)
\end{matrix} &
\begin{matrix}
(9/8)(4w-1)^2(12z+1)^2
\end{matrix}
\\
\noalign{\bigskip} 
\begin{matrix}
(1/2)\bigl[(4\pi^2+36\pi)z+(-8\pi^2+72\pi)w+\\
+\pi^2-15\pi+108\bigl]
\end{matrix} &
\begin{matrix}
(4w-1)^2(12\pi w-\pi+18)
\end{matrix}
\end{pmatrix}
$$
which has all entries $>0$ in the above range of values for $z$ and $w$.
\end{proof}

Of course the analogue of Lemma~\ref{ref7} (with $M'$ replacing $M$) holds for $q\in P_3\cup P_4$. We define a cone field $U_q$ on $T\oi^2$ as follows:
$$
U_q=\begin{cases}
(T_{M^{-1}q}M)[U^2_{M^{-1}q}], &
  \text{if }q\in MH\setminus(\text{closure of }NK); \\
(T_{{M'}^{-1}q}M')[U^2_{{M'}^{-1}q}], &
  \text{if }q\in M'H\setminus(\text{closure of }NK); \\
(T_{N^{-1}q}N)[U^4_{N^{-1}q}], &\text{if }q\in NK.
\end{cases}
$$
$U_q$ is defined $\lambda$-almost everywhere and depends piecewise-continuously on $q$. Let $\Sigma$ be the set of all $q$'s such that the cone field is defined and $B$ is differentiable on all the bilateral orbit of $q$: clearly $\Sigma$ has full measure. For $q\in\Sigma$ there are four possibilities:
\begin{itemize}
\item[(a)] $q,Bq\in (MH\cup M'H)\setminus NK$;
\item[(b)] $q\in (MH\cup M'H)\setminus NK$, $Bq\in NK$;
\item[(c)] $q\in NK$, $Bq\in (MH\cup M'H)\setminus NK$;
\item[(d)] $q,Bq\in NK$.
\end{itemize}
Using Lemma~\ref{ref7} one checks immediately that in all cases $(T_qB)[U_q]$ is contained in $U_{Bq}$, i.e., the cone field is invariant under the differential $T_qB$. Moreover, in all cases except~(a) the containment is strict. Almost every point $(r,\theta)\in H$ has an $F$-orbit which is dense in the circle of radius $r$. This implies that almost every $q\in\Sigma$ enters in $NK$ under $R_1^3$, and hence one of the cases~(b), (c), (d) must occur along the forward $B$-orbit of $q$. The cone field $U_q$ is therefore eventually strictly invariant, and this concludes the proof of Lemma~\ref{ref2}.

\section{Proof of Lemma~\ref{ref3}}

We use a geometric argument, analogous to those in~\cite[\S5]{BurtonEaston80} and \cite[Lemma~6]{wojtkowski80} for the toral case. Here the situation is more delicate, since we do not have the universal covering of the torus at our disposal. Nevertheless, we may use the homeomorphism $(r,\theta)\mapsto N(2r,\theta)$ from $D$ to $\oi^2$ to obtain a ramified covering of $\oi^2$ by $\oi\times\Rbb$ in the obvious way: keeping in mind this covering makes the argument below easier to follow.

The forward iterates of local unstable leaves are broken segments, and are all similar to each other (because all of them tend to approximate the unstable foliation). Analogously for the backwards iterates of stable leaves. The following pictures of such approximations (unstable and stable, respectively) are obtained for $l=m=1$, but clearly the pattern does not change for higher values of $l$ and $m$.

\pagebreak[4]

\begin{figure}[!h]
\begin{center}
\includegraphics[height=4cm]{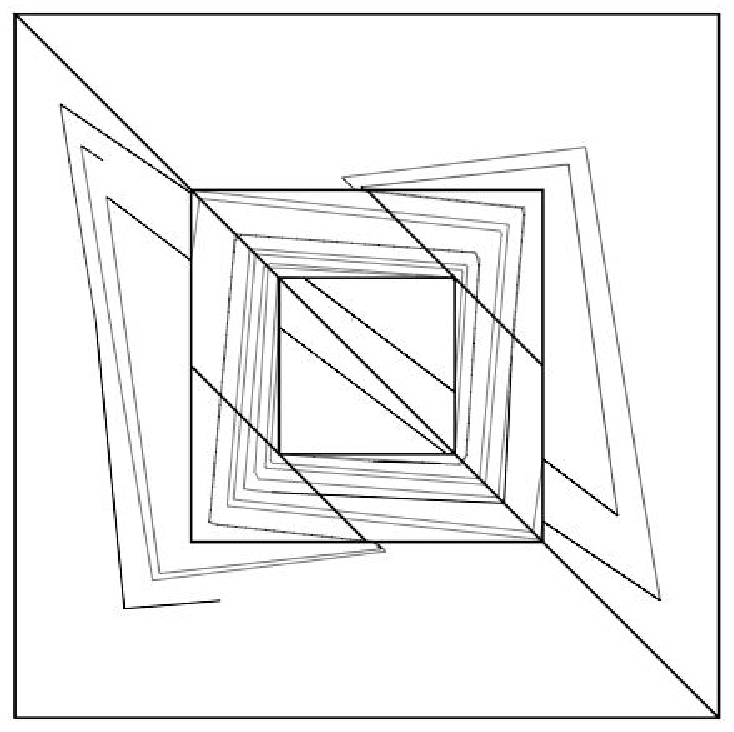}
\qquad
\includegraphics[height=4cm]{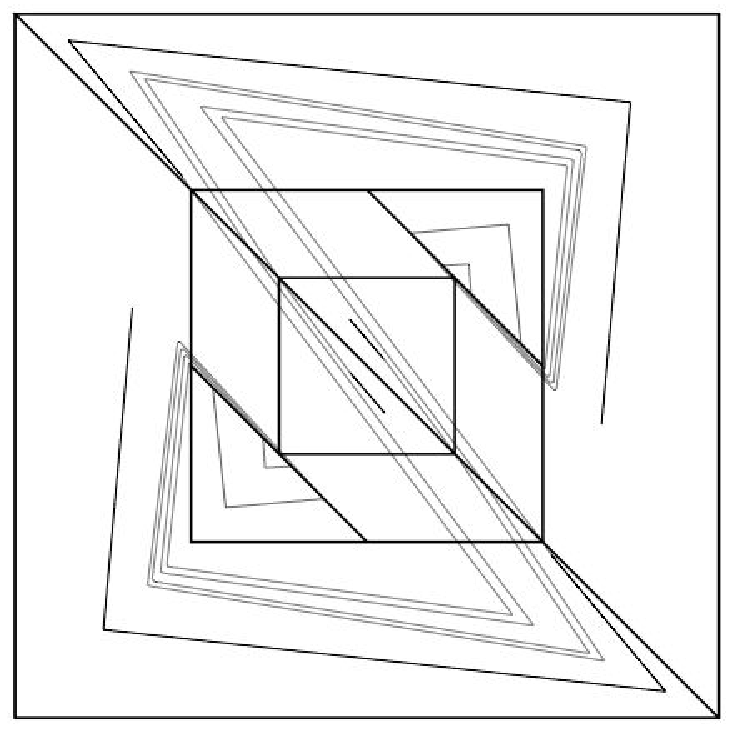}
\end{center}
\end{figure}

Let $P_i$ be one of the four parallelograms in the preceding section.
The \newword{horizontal sides} of $P_i$ are those fixed by $R_2$, and the \newword{vertical sides} those fixed by $R_1$.
A broken segment $\eta$ is \newword{horizontal in $P_i$} if some component of $\eta\cap P_i$ connects the vertical sides of $P_i$; we have a dual definition of $\eta$ being \newword{vertical in $P_i$}.
A simple geometric reasoning shows that if $\eta$ is horizontal in some $P_i$, then $B^k[\eta]$ is horizontal in all $P_i$'s, for every $k\ge2$ (if at least one of $l,m$ is $>1$, then $k\ge1$ is sufficient). Analogously, if $\eta$ is vertical in some $P_i$, then $B^{-h}[\eta]$ is vertical in all $P_i$'s for every $h\ge2$.

Since horizontal segments intersect vertical ones, it is sufficient to show that $B^k[\gamma^u]$ is horizontal in $P_i$, for some $i$ and $k\ge0$, and analogously $B^{-h}[\gamma^s]$ is vertical in some $P_i$; here $\gamma^u=\gamma^u_q$ and $\gamma^s=\gamma^s_p$ for $p,q\in\Sigma$ as in the preceding section. We treat only the unstable case, the stable one being similar. The length of $B^k[\gamma^u]$ goes to infinity as $k$ increases. Since the twists $R_1$ and $R_2$ ``go in the same direction'' (as in Case~1 of~\cite[Theorem~1]{wojtkowski80}), 
$B^k[\gamma^u]$ will not remain contained in any small region
(the dual phenomenon, in which the twists ``go in opposite directions'' and $B^k[\gamma^u]$ remains contained, is described in~\cite[p.~346]{przytycki83}). In particular, for $k$ sufficiently large, $B^k[\gamma^u]$ will cross some parallelogram, say $P_3$. The geometry of the unstable cone field forces $B^k[\gamma^u]$ to cross $P_3$ either horizontally, or through $\angles{q_3,q_2}$
and $\angles{p_4,q_3}$, or through $\angles{p_4,p'_1}$ and $\angles{p'_1,q_2}$. In the first case we are done, so suppose there is a component $\eta_0$ of $B^k[\gamma^u]\cap P_3$ that connects $\angles{q_3,q_2}$ with $\angles{p_4,q_3}$. Let $a_0$ be the extreme of $\eta_0$ that is in $\angles{p_4,q_3}$, and $b_0$ the extreme that is in $\angles{q_3,q_2}$. Compute $\eta_1=R_1^{3l}\eta_0$, $\eta_2=R_2^{4m}\eta_1=B\eta_0$, $\eta_3=R_1^{3l}\eta_2$, $\ldots$, and define:
\begin{itemize}
\item if $j$ is even, $a_{j+1}=a_j$, $b_{j+1}=$ the first point, moving from $a_{j+1}$ towards $R_1^{3l}b_j$ along $\eta_{j+1}$, in which $\eta_{j+1}$ crosses $\angles{q_3,q_2}$;
\item if $j$ is odd, $b_{j+1}=b_j$, $a_{j+1}=$ the first point, moving from $b_{j+1}$ towards $R_2^{4m}a_j$ along $\eta_{j+1}$, in which $\eta_{j+1}$ crosses $\angles{p_4,q_3}$.
\end{itemize}
The points $a_j$ and $b_j$ divide $\eta_j$ in three parts, a central one which remains in $P_3$ and approaches $q_3$, and two lateral ones, which get longer as $j$ increases. In particular the ``left part'' (the one that has $a_j$ as an extreme) will eventually cross $P_2$.
Fix now $j$ to be the first index for which this happens, and let $c_j$ be the first point, moving from $a_j$ along the left part, in which $\eta_j$ intersect a side of $P_2$. Again by the geometry of the unstable cone field, either $c_j\in\angles{p_2,q_0}$ and $j$ is even, or 
$c_j\in\angles{p_4,p_2}$ and $j$ is odd.
If $c_j\in\angles{p_2,q_0}$, then $\eta_j$ is horizontal in $P_2$ and we are done. If $c_j\in\angles{p_4,p_2}$, then the image under $R_2^{4m}$ of the part of $\eta_j$ between $c_j$ and $b_j$ is horizontal in $P_3$, and again we are done.

We left open the possibility that there is a component $\eta_0$ of $B^k[\gamma^u]\cap P_3$ that connects $\angles{p_4,p'_1}$ with $\angles{p'_1,q_2}$, but the proof is analogous. This time we take $a_0$ (respectively, $b_0$) to be the extreme of $\eta_0$ in $\angles{p'_1,q_2}$ (respectively, $\angles{p_4,p'_1}$), and construct $\eta_j,a_j,b_j$ as in the preceding case. The part of $\eta_j$ that has $a_j$ as an extreme will eventually cross the triangle $\angles{p'_0,p'_1,p'_2}$ and the parallelogram $P_4$. Again, $\eta_j$ must enter $P_4$ through $\angles{q_2,p'_2}$ and leave $P_4$ either at $c_j\in\angles{q_1,p_3}$ with $j$ even, or at $c_j\in\angles{q_2,q_1}$ with $j$ odd. But then either $\eta_j$ (in the first case) or 
the image under $R_2^{4m}$ of the part of $\eta_j$ between $b_j$ and $c_j$ (in the second case) crosses $P_4$ horizontally. This concludes the proof of Lemma~\ref{ref3} and hence of Theorem~\ref{ref1}.

\newcommand{\noopsort}[1]{}

\end{document}